\documentclass[11pt]{amsart}
\usepackage[left=1.0cm,top=2cm,right=2cm]{geometry}
\usepackage{bbold}
\usepackage[usenames,dvipsnames,svgnames,table]{xcolor}
\usepackage[pdftex,pagebackref,colorlinks=true,linkcolor=BrickRed,citecolor=OliveGreen,pdfstartview=FitH]{hyperref}
\usepackage{MnSymbol,wasysym}
\usepackage{graphicx}
\usepackage{ulem}

\newtheorem{thm}{Theorem}[section]

\newtheorem{thma}{Theorem}[section]

\newtheorem{coro}[thm]{Corollary}
\newtheorem{lem}[thm]{Lemma}
\newtheorem{prop}[thm]{Proposition}

\theoremstyle{definition}

\theoremstyle{remark}

\renewcommand{\index}{{\rm index}}

\newcommand{\oU}{\overline{U}}

\newcommand{\oz}{\overline{z}}

\newcommand{\ow}{\overline{w}}

\newcommand{\supp}{{\rm supp}}

\newcommand{\tu}{\widetilde{u}}

\newcommand{\tp}{\widetilde{p}}
\newcommand{\tC}{\widetilde{C}}

\newcommand{\Rbb}{ {\mathbb R}}

\newcommand{\Cbb}{ {\mathbb C}}

\newcommand{\Zcal}{ {\mathcal Z}}

\renewcommand{\phi}{\varphi}

\renewcommand{\th}{\theta}

\definecolor{darkolivegreen}{rgb}{0.33, 0.42, 0.18}

\newcommand{\cjin}[1]{\textcolor{blue}{#1}}

\title{Erratum: Euclidean Triangles have no hot spots.}
\author{Chris Judge and Sugata Mondal }
\date{June 2020}

\begin{document}

\maketitle


\section{Introduction}

We discovered a gap in the proof of the main theorem 
in our paper `Euclidean triangles have no hot spots' \cite{Annals}.  
(\cjin{ \href{https://arxiv.org/abs/1802.01800v4}{Version 4} 
of this arxiv posting is equivalent to the published version.})
In particular, Lemma 3.4 in that paper is false as stated, and hence
Lemma 7.6 is unproven. Though we strongly believe that the statement 
of Lemma 7.6 is true, we do not currently know of a proof.
The proof of Theorem 1.1 \cite{Annals} is still valid for obtuse triangles
since a weaker version of Lemma 7.6 suffices 
(see \S \ref{sec:obtuse-right} below). 
The proof of Theorem 1.1 \cite{Annals}
is however incomplete in the case of acute triangles. 
Nonetheless, using an alternative argument based on \cite{Annals},
we prove here that acute triangles have no hot spots:

\begin{thma}\label{thm:main}
If $u$ is a second Neumann eigenfunction of the Laplacian
on an acute triangle $T$,
then $u$ has no critical points in the interior of $T$. 
\end{thma}

The bulk of this erratum consists of a proof of Theorem A
using the same general approach taken in \cite{Annals}. 
In the last section, \S \ref{sec:obtuse-right}, we explain the minor modifications
of \cite{Annals} needed to recover Theorem 1.1 of \cite{Annals} in the 
case of obtuse and right triangles. The results in \S 13 of \cite{Annals} 
are also correct in the case of obtuse or acute triangles.
So, in fact, by combining the results in Theorem 1.1, 
Theorem 13.4, and, Corollary 13.5 in \cite{Annals}, we prove 
that each second Neumann eigenfunction of an obtuse or right triangle 
has no critical points on $T$ 
(see Theorem \ref{thm:obtuse-right} below).

In addition to Theorem 1.1, Lemma 3.4, and 
Lemma 7.6  in \cite{Annals}, the proofs of the following statements 
in \cite{Annals} are incomplete: Theorem 7.7, 
Corollary 7.8, Proposition 10.5, Lemma 11.2, Theorem 13.1, Corollary 13.2,
Proposition 13.3\footnote{The statement of Proposition 13.1 has been proven 
in \cite{Miyamoto} but the proof in \cite{Annals} is incomplete.}
Theorem 13.4,
and Corollary 13.5.
(As indicated above, the proofs of the results of \S 13 
are valid for right and obtuse triangles.)
Except for Lemma 3.4, each of these propositions 
would follow from Lemma 7.6 which we believe to be true.



\section{Preliminaries}

As in \cite{Annals}, the symbol $\mu$ 
denotes the second Neumann eigenvalue of a triangle $T$, 
and $u$ denotes an associated second Neumann 
eigenfunction. Also,
we suppose that the triangle $T$ lies in the complex plane, 
and we let $z=x+iy$ denote  a point in $\Cbb \cong \Rbb^2$.

\begin{prop}  \label{lem:pos-local-max}
Let $u$ be a second Neumann eigenfunction of a polygon $P$.  
An acute vertex $v$ of $P$ is a local extremum of $u$ 
if and only if $u(v) \neq 0$.
Moreover, a vertex $v$ is a local maximum (resp. minimum) of $u$
if and only if $u(v)>0$ (resp. $u(v)<0$).  
\end{prop}

\begin{proof}
Since the angle at $v$ is acute, 
the first statement follows from  Lemma 4.3 and Corollary 5.3 \cite{Annals}.
In particular, if $u$ is a local extremum, then from the discussion in 
\S 4 \cite{Annals}, we have that
$u(z) = u(v) \cdot J_0 \left(\sqrt{\mu} \cdot |z-v| \right) + o(|z-v|^2)$
where here $J_{\nu}$ denotes the Bessel function of the first kind and of order $\nu$.
It is well known that $J_0(0)=1$, $J_0'(0)=0$, and  $J_0''(0)<0$.
The claim follows.
\end{proof}

As in \cite{Annals}, we let $\Zcal(f)$ denote the set of $z$ such that $f(z)=0$,
and we let $L_e$ denote a constant vector field that is parallel to the side $e$.

By Lemma 6.1 and Proposition 6.2 in \cite{Annals}, the nodal set $\Zcal(L_e u)$ may be naturally regarded 
as a locally finite  tree whose vertices consist of critical points of $u$ and, possibly, some vertices of $T$.
The {\it degree of a vertex} $p$ in the graph $\Zcal(L_e u)$ is the number 
of edges (arcs) that have $p$ as an endpoint.

\begin{lem}  \label{lem:vanish-vertex-nodal}
Let $v$ be a vertex of an acute triangle $T$, let $e$ be the side  opposite to $v$, and
let $e'$ be a side that is adjacent to $v$.
If $u(v)=0$, then
\begin{enumerate} 
\item  $v$ is not a degree 1 vertex of $\Zcal(L_e u)$, and 
\item  $v$ is a degree 1 vertex of $\Zcal(L_{e'}u)$.
\end{enumerate}

If $u(v) \neq 0$, then
\begin{enumerate} 
\item  $v$ is a degree 1 vertex of $\Zcal(L_e u)$, and 
\item  $v$ is not a degree 1 vertex of $\Zcal(L_{e'}u)$.
\end{enumerate}
\end{lem}

\begin{proof}
By applying a rigid motion to $T$ if necessary, we may assume, without loss of generality,
that $v=0$, the adjacent side $e'$ lies in the positive real axis, and $e$ 
lies in the upper half plane.
We will use the Bessel expansion of $u$ 
as described in equation (6) of \cite{Annals}. 

If $u(v)=0$, then using Lemma 4.3 and Corollary 5.3 \cite{Annals}, we find that
\[ 
u(z)~ 
=~
a \cdot |z|^{\nu} \cdot \cos(\nu \cdot \arg(z))~
+~
o(|z|^\nu)
\]
where $a \neq 0$ and $\nu = \pi/\beta>2$ where $\beta$ is the angle at $v$.
Without loss of generality, the vector field 
$L_e$ has unit norm, and hence there exists $\psi \in \Rbb$
so that $L_e= \cos(\psi) \cdot \partial_x + \sin(\psi) \cdot \partial_y$.
Because $\partial_x = \cos(\theta) \cdot \partial_r - \sin(\th) r^{-1} \cdot \partial_{\th}$ 
and $\partial_y = \sin(\theta) \cdot \partial_r + \cos(\th) r^{-1} \cdot \partial_{\th}$,
we have 
\begin{equation} \label{eqn:Lpsi-polar}
L_e~
=~ \cos(\psi-\theta) \cdot \partial_r~
+~
\sin(\psi-\theta) \cdot \frac{1}{r} \partial_{\theta}.
\end{equation}
A straightforward computation gives
\begin{equation} \label{eqn:L-u}
L_{e} u \left(r e^{i \th} \right)~
=~
a \cdot \nu \cdot r^{\nu-1} \cdot \cos(\psi - \th +\nu \cdot \th)~
+~
o(r^{\nu-1}),
\end{equation}
where $r = |z|$ and $\theta = \arg(z)$.

Because $T$ is a acute and $\psi$ is the exterior angle at the vertex 
where $e$ and $e'$ meet, we have $\pi/2 <\psi< \pi/2 + \beta$, and hence
$\pi/2 < \psi - \th +\nu \cdot \th < 3\pi/2$ for $\theta \in [0, \beta]$.
Thus $L_e u$ does not vanish in a neighborhood of $v$, and in particular
$v$ is not a degree 1 vertex of $\Zcal(L_e u)$. 

Similarly, we may assume that $L_{e'} = \partial_x$,
and so a straightforward computation gives
\begin{equation*}
  L_{e'} u \left(r e^{i \th} \right)~
=~
a \cdot \nu \cdot r^{\nu-1} \cdot \cos((1 -\nu) \cdot \th)~
+~
o(r^{\nu-1}).  
\end{equation*}
Since $\nu > 2$ it follows that $L_{e'} u$ changes sign as $\th$ varies over $[0, \beta]$.
In particular, $v$ is a vertex of $\Zcal(L_{e'} u)$. 
Another explicit computation for $\partial_\th L_{e'} u (r, \th)$ shows that $\partial_\th L_{e'} u (r, \th)$ does not vanish for $r$ sufficiently small. This provides that $v$ is a degree one vertex of $\Zcal(L_{e'} u)$.

For the second part, if $u(v) \neq 0$, then from equation (6) of \cite{Annals} we find that
\[ 
u(z)~ 
=~ a~ -~ 
b \cdot |z|^2 ~
+~
o(|z|^2)
\]
where $a \neq 0 \neq b$.
Computing as above we get
\begin{equation}
    L_e u \left(re^{i\theta}\right)~ =~ - b \cdot \cos(\psi-\theta) \cdot r~ 
    +~  o(r).
\end{equation}
Since $\pi/2 < \psi < \pi$ we get the first assertion of the second part.

For the second part of the second assertion we may again assume that $L_{e'} = \partial_x$. Hence we have
\begin{equation}
    L_{e'} u \left(re^{i\theta} \right)~ =~  - b \cdot \cos(\theta) \cdot r~ +~ o(r).
\end{equation}
since $\theta$ ranges on the interval $[0, \beta]$ and $\beta < \pi/2$, we get the second assertion of the second part.
\end{proof}

The following is, in some sense, a corollary of Lemma 7.5 \cite{Annals}.

\begin{prop} \label{prop:degree-two-implies-degree-zero}
Let $e$, $e'$ and $e''$ be the sides of an acute triangle $T$, and let 
$c$ be a critical point of $u$ that lies on $e$.
Suppose that $c$ is a degree two vertex of $\Zcal(L_e u)$. 
Then $c$ is degenerate, and, 
if, as a vertex of $\Zcal(L_{e'} u)$, the degree of $c$ equals two,
then, as a vertex of $\Zcal(L_{e''} u)$, the degree of $c$ equals zero.

\end{prop}

\begin{proof}
As in \S 7 \cite{Annals}, we may assume that $e$ lies in the 
real axis and $c=0$. Moreover, since $u$ is a Neumann function
it extends by reflection to an eigenfunction in a neighborhood of $c$.
By Lemma 7.2 \cite{Annals}, we have\footnote{
Note that the inclusive `or' is intended in the last sentence of the statement of Lemma 7.2 \cite{Annals}.}
$u(z)= a_{00} + a_{20} \cdot x^2 + a_{02} \cdot y^2 + o(|z|^2)$
and hence $\partial_x u(z)= 2 a_{20} \cdot x + o(|z|)$.
Thus, since $c$ is a degree $2$ vertex of $\Zcal(L_e u)$, we have 
$a_{20}=0$ and in particular $c$ is a degenerate critical point.

Since $c$ is a degree two vertex of $\Zcal(L_e u)$, it follows that the Taylor expansion of $u$ 
at $c$ has the form 
\[ 
u(z)~
=~
a_{00}~ +~ a_{02} \cdot y^2~ +~ a_{30} \cdot (x^3 - 3 x \cdot y^2) +  o(|z|^3) \] 
where $a_{02} \neq 0 \neq a_{30}$.
We have $L_{e'} = c_1' \partial_x + c_2' \partial_y$
and $L_{e''}=c_1'' \partial_x + c_2'' \partial_y$
for constants $c_1', c_2', c_1'', c_2''$.  Since the triangle
$T$ is acute, none of the constants equal zero, and, we may assume that 
the constants $c_1'$ and $c_1''$ have opposite signs, and 
the constants $c_2'$ and $c_2''$ have the same sign.
By Lemma 7.5 \cite{Annals}, since $c$ is a degree two vertex of $L_{e'} u$
we have $(a_{30}/a_{02}) \cdot (c'_1/c'_2)<0$. 
Hence $(a_{30}/a_{02}) \cdot (c''_1/c''_2)>0$, and so, by 
Lemma 7.5 \cite{Annals}, the point $c$ is a degree 0 vertex of 
$\Zcal(L_{e''} u)$. 

\end{proof}


\subsection{The index of a critical point}

We adapt the notion of the `index' of a critical point
in the sense of the Poincar\'e-Hopf index theorem.
To do this, we let $DT$ be the double of the triangle $T$, 
the conically singular surface obtained by gluing 
two isometric copies of $T$ along the boundary using the identity map.
The Neumann eigenfunction $u$ lifts to a function $\tu: DT \to \Rbb$.
For each point $\tp \in DT$, the level set of $\tu$ that contains $\tp$
is a finite graph $G_{\tp}$. If $p \in T$ corresponds to $\tp$, then define 
the {\it index of $p$} to be $1-k/2$ where $k$ is the valence of 
$G_{\tp}$  at $\tp$. If $\nabla u(p)=0$,
then the index of $p$ equals the Poincar\'e-Hopf 
index of $\nabla u$ at its zero $p$.

The lift $\tu$ is smooth on the complement of the conical singularities
corresponding to the vertices of $T$. 
It follows that if the index of $p$ is nonzero, then $p$ is either a vertex of $T$
or $p$ is a critical point of $u$. Note that $p$ is a local extremum
if and only if the index of $p$ equals $+1$. If the index of a critical
point equals $-1$, then we will call $p$ a {\it saddle}. 

Let $C^{\partial}$ denote the set of $p \in \partial T$ whose index is nonzero,
and let $C^{\circ}$ denote the set of $p \in T^\circ$ whose 
index is nonzero. Corollary 5.7 \cite{Annals} implies that each of 
these sets is finite. 

\begin{prop}[Index formula]  \label{prop:index-thm}
Let $u: T \to \Rbb$ be a Neumann eigenfunction. Then
\[
2~
=~ 
\sum_{p \in C^{\circ}}
2 \cdot \index(p)~
+~
\sum_{p\, \in\, C^{\partial}}
\index(p).
\]
\end{prop}

\begin{proof}
For each $\tp \in DT$, define the index of $\tp$ to be the index of $p$.
A point $p \in \partial T$ corresponds to one point in $DT$,
and a point $p \in T^{\circ}$ corresponds to two points in $DT$.
Hence, it suffices to show that $\sum_{C} \index(\tp) =2$ where $\tC$ 
is the set of $\tp \in DT$ that have nonzero index.

Let $A$ be the union of the graphs $G_{\tp}$ over $\tp \in \tC$. 
The complement of $A$ consists of topological annuli, 
and hence, by Euler-Poincar\'e formula, $\chi(DT)=\chi(A)$. 
On the other hand, the number of edges in $A$
equals $\frac{1}{2} \sum_{\tp \in \tC} k_{\tp}$ where $k_{\tp}$ 
is the valence of the graph $A$ at $\tp$. It follows that 
$\chi(DT)= \sum_{\tp \in \tC} \index(p)$. Since $DT$ is a 
topological sphere, we have $\chi(DT)=2$, and the formula follows. 
\end{proof}

\begin{prop} \label{prop:index-at-least-minus-one}
The index of a critical point $p$ of a second 
Neumann eigenfunction $u$
is either $1,0$, or $-1$.
\end{prop}

\begin{proof}
Suppose not. Let $\tu$ be the lift of $u$ to the surface DT,
and let $\tp \in DT$ correspond to $p$.
Then the graph $\tu^{-1}(\tu(\tp))$ has degree $k > 4$ at $\tp$. 
It follows that, in the natural coordinates at $\tp$, 
we have $\tu(z)-\tu(\tp)= o(|z-\tp|^2)$.
In particular, the degree two homogeneous polynomial $h_2$ 
consisting  of second order terms in the Taylor expansion of $\tu$
at $p$ vanishes indentically. But $(\Delta h_2)(\tp) = \mu \cdot \tu(\tp)$
and so $u(\tp)=0$. 
But by Theorem  5.2 \cite{Annals}, the nodal set $\tu^{-1}(0)$ 
is a simple closed curve, and hence the degree $k$ equals two,  
a contradiction.
\end{proof}

As in \cite{Annals}, we will often consider the Taylor expansion of the eigenfunction $u$
in a neighborhood of a critical point that lies on a side $e$ of $T$. 
By applying a rigid motion if necessary, 
we will always assume that $e$ belongs to the real axis and that $p=0$.
Because $u$ is a Neumann eigenfunction, 
the reflection principle provides an extension $u$ to a neighborhood
of $0$ so that $u(\oz)= u(z)$. This extension will be implicit in
what follows.

\begin{lem} \label{lem:tangential-cusp}
Suppose that $p$ is an index zero critical point that belongs to the side $e$. 
Then there exist real-analytic functions $c: \Cbb \to \Rbb$ and $\rho:\Rbb \to \Rbb$ 
and an odd integer $k \geq 3$ so that $c(0) \neq 0$, $\rho(0)\neq 0$, and
\begin{equation} \label{eqn:tangential-cusp-equation} 
u(z)~
=~
u(0)~
+~
c(z) \cdot 
\left(  
y^2~
- x^k \cdot \rho(x)
\right).
\end{equation} 

\end{lem}

\begin{proof}
Because the index of the critical point $p$ of $u$ equals zero, 
the Hessian of $u$ has exactly one nonzero eigenvalue.
The eigenspace $E$ that corresponds to the nonzero eigenvalue is 
invariant under the reflection $z \mapsto \oz$. 
Thus $E$ is either the real or imaginary axis.

We claim that $E$ is not the real axis. 
Indeed, suppose to the contrary that $E$ is the real axis. 
Then $\partial_x u(0)=0$ 
but $\partial_x^2\, u(0) \neq 0$. 
The Weierstrass preparation theorem
applies to provide unique real-analytic functions $a$, $b_1$, and $b_2$
defined near $0$, so that $a(0) \neq 0$, $b_1(0)=0=b_2(0)$, and 
\[
u(z)~ -~ u(0)~ =~ a(z) \cdot \left(x^2~ +~ b_1(y)\cdot x~ +~ b_2(y) \right)
\]
for $z$ near $p=0$. Since the factorization is unique and $u(\oz)=u(z)$,
we have $b_j(y)= b_j(-y)$ for $j=1,2$. In particular, the 
discriminant $D(y):= b_1(y)^2 - 4 \cdot b_2(y)$ is an even function. 
If $D$ were to vanish on a neighborhood of $0$, then a straightforward computation 
shows that $\nabla u$ vanishes along the level set of $u$ that contains $p=0$. 
But the set of critical points of $u$ is finite by Corollary 5.7 \cite{Annals},
and so $0$ is an isolated zero of $D$. 
Since the critical point $p$ of $u$ has index zero  
and $\partial_x^2 u(0) \neq 0$, the level set $u^{-1}(u(p))$ near $p=0$
is an arc that meets either the upper half plane or lower plane. 
Since $D$ is even, it follows that $D(y) >0$ for $y \neq 0$ sufficiently small, 
and hence there exists a neighborhood $U$ of $p=0$ so that 
the intersection of $u^{-1}(u(p))-\{p\}$ and $U$ 
consists of four arcs. This contradicts the assumption that $p$ is a 
zero index critical point of $u$.

Therefore $E$ coincides with the imaginary axis. By use of 
the Weierstrass preparation theorem we find that 
\begin{equation} \label{eqn:weierstrass-cusp}
u(z)~ -~ u(0)~ =~ a(z) \cdot \left(y^2~ +~ b_1(x)\cdot y~ +~ b_2(x) \right)
\end{equation}
for unique real-analytic functions $a$, $b_1$, and $b_2$
defined near $0$ where $a(0) \neq 0$ and $b_1(0)=0=b_2(0)$.
Since the factorization is unique and  $u(\oz)=u(z)$, we find that
$b_1(x)= 0$. Because $u$ is an isolated critical point, there exists 
$\epsilon>0$ so that $b_2(x) \neq 0$  if $0 < |x|< \epsilon$.
We claim that, moreover, $b_2(x) \cdot b_2(-x)<0$ if $0 < |x|< \epsilon$.
we have $b_2(x) \cdot b_2(-x)<0$. Indeed, otherwise $b_2(x) \cdot b_2(-x)>0$, 
and thus from (\ref{eqn:weierstrass-cusp}) we find that  
there exists a neighborhood $U$ of $p=0$ so that 
the intersection of $u^{-1}(u(p))-\{p\}$ and $U$ 
consists of four arcs. This contradicts the assumption that $p$ is a 
zero index critical point of $u$.

Since $b_2(x) \cdot b_2(-x)<0$ the first nonzero term in the Taylor 
series of $b_2$ about zero has odd degree $k$, and since $\partial_x u(0)=0$
we also have $k \geq 3$. The claim follows. 
\end{proof}

\begin{prop} \label{prop:degree-1-stable}
Let $X$ be a constant vector field. If $p$ is a degree 1 vertex of 
$\Zcal(X u)$ that is not a vertex of $P$, then $p$ is 
a critical point with nonzero index. 
\end{prop}

We remark that a similar proof shows that
Proposition \ref{prop:degree-1-stable} holds when 
$X$ is a rotational vector field, but we will 
not use this fact in the erratum.

\begin{proof}
Since $p$ is not a vertex, $p$ lies in the interior 
of a side $e$. If the vector $X(p)$ were 
orthogonal to $e$, then $Xu$ would vanish on $e$,
and hence $p$ would not be a degree 1 vertex of $\Zcal(Xu)$.
Therefore $X(p)$ is independent of the outward
normal vector at $p$,
and in particular $p$ is a critical point of $u$. 
It remains to show that $p$ has nonzero index.

As above, we may suppose without
loss of generality that $e$ lies in the real-axis and that $p=0$.
Thus, since $X$ is not orthogonal to $e$, we may assume
without loss of generality that $X= \cos(\psi) \partial_x + \sin(\psi) \partial_y$ 
where $\psi \neq \pi/2 \mod \pi$.  

Suppose to the contrary that the index of $p$ were to equal zero.
Then by Lemma \ref{lem:tangential-cusp}, near $p$, the function
$u$ would satisfy (\ref{eqn:tangential-cusp-equation}) 
where $c(0) \neq 0 \neq \rho(0)$ and $k \geq 3$ is odd. 
Direct computation 
shows that $\partial_y X u(0)= 2 \cos(\psi) \cdot c(0)$
and hence $\partial_y X u(0) \neq 0$.  Thus, by the implicit function
theorem, there exists a function $f:(- \epsilon, \epsilon) \to \Rbb$
so that $f(0)=0$ and 
\begin{equation} \label{eqn:implicit-arc}
Xu \left(x + i \cdot f(x) \right)~ =~ 0. 
\end{equation}
From (\ref{eqn:tangential-cusp-equation}), we find that for each real $x$
\begin{equation} \label{eqn:cusp-real-expansion}
 (Xu)(x~ +~ i \cdot f(x))~ 
 =~ 
 \sin(\psi) \cdot c(0) \cdot 2 f(x)~
 -~
 \cos(\psi) \cdot c(0) \cdot k \cdot x^{k-1}~ 
 +~ 
 o\left(|x|^{k-1}\right)~
 +~
o\left(|f(x)| \right). 
\end{equation}
Therefore, since $\psi \neq \pi/2 \mod{\pi}$, we find that 
$f(x) = \frac{k}{2} \cdot \cot(\psi) \cdot x^{k-1} + O(|x|^k)$.
Since $k-1$ is even and greater than 0, 
the function $f$ is either positive in a deleted neighborhood 
of $0$ or negative in a deleted neighborhood of $0$. 
Thus, there exists a neighborhood $U$ of $p=0$ such that 
either $(U \cap \Zcal(X u))\setminus \{p\}$
lies in the upper half plane or  
$(U \cap \Zcal(X u))\setminus \{p\}$
lies in the the lower half plane. 
Hence $p$ is not a degree 1 vertex, a contradiction.
\end{proof}

Let $e$ be a side of $T$.  
Lemma 6.6 \cite{Annals} and Lemma \ref{lem:vanish-vertex-nodal}
imply that the set $e \cap \Zcal(L_e)$ consists 
of the critical points of $u$ and the endpoints of $e$ at which $u$ vanishes.

\begin{coro} \label{coro:n-arcs-n-1-stable}
Let $s$ be the number of nonzero index 
critical points on $e$, and let $t$ be the number of index zero
critical points on $e$.  Then $\Zcal(L_e u) \setminus e$ has at least 
$s + 2t$ components whose closures intersect $e$,
and the number of degree 1 vertices of 
$\Zcal(L_e u)$ that lie in $\partial T \setminus e$
is at least $s+2t$. 
The intersection of the closure 
of each component and $e$ consists of a critical point of $u$.

If $u$ vanishes at an endpoint $v$ of $e$,
then $\Zcal(L_e u) \setminus e$ has at least 
one additional component with at least one
degree 1 vertex in $\partial T \setminus e$, and thus there exist
at least $s + 2t+1$ components whose closures intersect $e$.
Moreover, the intersection of the closure of the component and $e$
consists of the point $v$.
\end{coro}

\begin{proof}
By Lemma 6.6 \cite{Annals}, each critical point $p$ in $e$ is the endpoint of 
at least one arc in $\Zcal(L_e u)$ that intersects the interior of $T$. 
If the index of $p$ is zero then, it follows from Lemma \ref{lem:tangential-cusp} 
that $p$ is the endpoint of at least two arcs in $\Zcal(L_e u)$ that intersect
the interior of $T$. If $u$ vanishes at a vertex  adjacent to $e$, then 
Lemma \ref{lem:vanish-vertex-nodal} implies that this vertex is a
degree 1 vertex of  $\Zcal(L_e u)$.

Lemma 3.3 \cite{Annals} implies that the closure of each component 
of  $\Zcal(L_e u) \setminus e$ intersects $e$ in at most one point. 
Since a finite tree has at least two degree 1 vertices, 
each component $\Zcal(L_e u) \setminus e$ has at least one vertex in 
$\partial T \setminus e$.
\end{proof}

\begin{prop} \label{prop:crit_opposite_sign_on_vertices}
Let $v'$ and $v''$ be the vertices of $T$ adjacent to the side $e$.
If $u(v') \cdot u(v'') <0$, then the number of nonzero index critical 
points lying in the interior of $e$ is even.   
If $u(v') \cdot u(v'') >0$, then the number of nonzero index critical 
points lying in the interior of $e$ is odd.   
\end{prop}

\begin{proof}
Let $\gamma:[0,1] \to e$ be a parameterization of $e$ such that 
$\gamma(0)=v'$ and $\gamma(1)=v''$. 
Since $u$ satisfies Neumann conditions, 
a point $p:=\gamma(t)$ is a critical point of $u$ if and only $t$ it is a critical 
point of $u \circ \gamma$. 

We claim that $t$ is a local extremum of $u \circ \gamma$
if and only if the index of the critical point $p=\gamma(t)$ of $u$ is nonzero. 
Indeed, by Proposition \ref{prop:index-at-least-minus-one},
the index of $u$ at $p =\gamma(t)$ is either $+1$, $0$, or $-1$.
If the index of $u$ at $p$ equals $+1$, then $p$ is a local 
extremum of $u$, and hence $t$ is a local extremum of $u \circ \gamma$.
If the index of $p$ equals zero, then, by Lemma \ref{lem:tangential-cusp}, it follows
that $t$ is not a local extremum.
If the index of $\gamma(t)$ equals $-1$, then there is a neighborhood $U \subset T$
of $p = \gamma(t)$ and $U$ contains only one critical point of $u$.
It follows that $t$ is a local extremum of $u \circ \gamma$ if the 
index is $-1$.

If $u(v')$ and  $u(v'')$ have opposite signs, then  
by replacing $u$ by $-u$ if necessary, we may assume that $u(v') < 0 < u(v'')$. 
Proposition \ref{lem:pos-local-max} then implies that the composition $u \circ \gamma$ 
is increasing near both $0$ and $1$.  Hence the number of local extrema of 
$u \circ \gamma$ is even, and thus the number of nonzero index critical points of $u$
that lie on $e$ is also even.

The proof of the other assertion is similar.
\end{proof}

The following is a strengthening of Lemma 7.1 \cite{Annals}.

\begin{prop} \label{prop:at-least-four}
If an eigenfunction $u$ has an interior critical point, then 
$u$ has at least four critical points with nonzero index
that lie on the boundary of $T$.
\end{prop}

\begin{proof}
Since $u$ has an interior critical point, 
Lemma 7.1 \cite{Annals} implies that $\Zcal(L_e u)$ has at least one degree 1 vertex 
in the interior of each side of $T$. By Proposition \ref{prop:degree-1-stable}, each of these degree one vertex is a critical point of $u$ with nonzero index. 
By Theorem 5.2 \cite{Annals}, the nodal set $\Zcal(u)$ is a simple arc 
with two endpoints in $\partial T$ that divides $T$ into two components. 
By Corollary 5.5 \cite{Annals} at most one endpoint of $\Zcal(u)$ 
is a vertex of $T$. Let $e$ be a side of $T$ whose interior 
contains an endpoint of $\Zcal(u)$.  By Lemma 5.1 \cite{Annals}, 
the function $u$ does not vanish at either of the vertices, $v'$, $v''$,
that are adjacent to $e$.  In particular, the vertices $v'$ and $v''$ lie 
in distinct components of $T \setminus \Zcal(u)$, and hence
$u(v')$ and $u(v'')$ have opposite signs. 
Since $e$ already contains one nonzero index critical point, 
Proposition \ref{prop:crit_opposite_sign_on_vertices}
implies that the side $e$ contains at least two nonzero index critical points.
Thus, we have at least four nonzero index critical points in total. 
\end{proof}


\subsection{The boundary integral}

Integration by parts implies that the second Neumann eigenfunction $u$ satisfies 
$$
\int_T |\nabla u|^2~ 
=~ 
\mu \int_T |u|^2.
$$ 
However, if, for example, $L$ is a constant vector field, then the 
eigenfunction $Lu$ 
does not satisfy a Neumann condition, and so 
$$
\int_T |\nabla Lu|^2~
=~
\mu \int_T |Lu|^2 + \int_{\partial T} Lu \cdot \partial_{\nu}Lu
$$
where $\partial_{\nu}$ is the outward normal derivative. We collect here
some results that will help us determine the sign of the `boundary integral'
$\int_{\partial T} Lu \cdot \partial_{\nu}Lu$.

The following is a variant of Proposition 3.1 in \cite{L-R}. 

\begin{lem} \label{lem:boundary-integral-positive}
Let $U \subset T$ such that $\partial U \cap \partial T$ 
lies in a line $\ell$ and $T \setminus U$ has nonempty interior. 
Let $\phi \in H^1(T)$, $\supp(\phi) \subset \oU$,
$\phi|_{U} \in H^2(U)$, and  
$\Delta \phi|_U= \mu \cdot \phi|_U$.
Then  
\[ 
\int_{\partial T} \phi \cdot \partial_{\nu} \phi~ 
>~ 
0.
\]
\end{lem}

\begin{proof}
Suppose to the contrary that the boundary integral of $\phi$ is nonpositive.
Then since $\phi$ vanishes on $\partial U \cap T^{\circ}$
and $\Delta \phi|_U = \mu \cdot \phi|_U$,
integration by parts gives
\[
\int_{T} |\nabla \phi|^2~ 
=~
\int_{U} \nabla \phi \cdot \nabla \phi~  
=~
\mu \int_{U} |\phi|^2~
+~
\int_{\partial U \cap \partial T} \phi \cdot \partial_{\nu} \phi~ 
\leq~
\mu \int_{T} |\phi|^2.  
\]

Without loss of generality $\partial U \cap \partial T$ lies 
in the real-axis. If we let
$w(z)=\exp \left(i \cdot \sqrt{\mu} \cdot x \right)$,
then $\nabla w = i \sqrt{\mu} \cdot w \cdot \partial_x$ 
and hence $\int |\nabla w|^2 = \mu \int |w|^2$.
Moreover, along the real axis $\partial_{\nu} w = 0$  
and so integration by parts gives
\[
\int_{T} \nabla \phi \cdot \nabla \ow~
=~
\int_{U} \nabla \phi \cdot \nabla  \ow~ 
=~
\int_{U} \phi \cdot  \Delta \ow~
+~
\int_{\partial U} \phi \cdot \partial_{\nu} \ow~
=~
\mu
\cdot
 \int_{T} \phi \cdot \ow.
\]
By combining the above, we find that for each pair  $a$, $b \in \Cbb$
\[
\int_{T} |\nabla (a \cdot \phi~ +~ b \cdot w) |^2~
\leq~  
\mu \cdot \int_{T} |a \cdot \varphi~ +~ b \cdot w|^2.
\]
Since $\phi$ vanishes on an open set, the functions $\phi$ and $w$
are linearly independent. Thus, since $\mu$ is the second Neumann 
eigenvalue, there exist $a$ and $b$ such that 
$a \cdot \phi + b \cdot w$ is a Neumann eigenfunction associated to $\mu$.
The function $w$ does not satisfy Neumann conditions on the entire
boundary of $T$ and so $a \neq 0$. Thus, since $\phi$ vanishes on an open set,
the function $a \cdot \phi + b \cdot w$ is not real-analytic on $T$, 
a contradiction.
\end{proof}

The following lemma expresses an idea found in \S 3 of \cite{Miyamoto}.

\begin{lem}  \label{lem:boundary-term-nonneg}
Let $\Omega$ be a bounded plane domain with Lipschitz boundary,
and let $U_-, U_+ \subset \Omega$ be two disjoint nonempty
open subsets with Lipschitz boundary.
Let $\phi_-, \phi_+ \in H^1(\Omega)$ be two real-valued functions 
with $\supp(\phi_{\pm}) \subset \oU_{\pm}$.
If the restrictions $\phi_{\pm}|_{U_{\pm}}$ are smooth and satisfy 
$\Delta \phi_{\pm}= \mu \cdot \phi_{\pm}$, then the integrals 
\begin{equation} \label{eqn:boundary-integrals}
\int_{\partial \Omega} \phi_{+} \cdot \partial_{\nu} \phi_{+}
\mbox{ \   and   \   } 
\int_{\partial \Omega} \phi_{-} \cdot \partial_{\nu} \phi_{-}
\end{equation}
cannot both be negative. 
\end{lem}

\begin{proof}
Since $\Delta \phi_{\pm} =\mu \cdot \phi_{\pm}$ 
and $\supp(\phi) \subset \oU$, integration by parts gives
\begin{equation}  \label{eqn:phi-energy}
\int_{\Omega} \nabla \phi_{\pm} \cdot \nabla \phi_{\pm}~
=~
\int_{U_{\pm}} \phi_{\pm} \cdot \Delta \phi_{\pm}~ 
+~
\int_{\partial U_{\pm}} \phi_{\pm} \cdot \partial_{\nu} \phi_{\pm}~
=~
\mu \cdot  \int_{U_{\pm}} |\phi_{\pm}|^2~
+~
\int_{\partial \Omega_{\pm}} \phi_{\pm} \cdot \partial_{\nu} \phi_{\pm}.  
\end{equation}
Since $\phi_+$ and $\phi_-$ have disjoint support, 
$\int \nabla \phi_{+} \cdot \nabla \phi_- = 
0 = \int \phi_{+} \cdot \phi_-$. Hence for each $a_+,a_- \in \Rbb$ 
\[
\int_{\Omega} \left| a_+ \cdot \nabla \phi_+ + a_- \cdot \nabla \phi_-\right|^2~
=~
\mu_2 \cdot \int_{\Omega} |a_+ \cdot \varphi_+~ +~ a_- \cdot \phi_-|^2~
+~
a_+^2 \cdot \int_{\partial \Omega} \phi_+ \cdot \partial_{\nu} \phi_+~
+~ 
a_-^2 \cdot \int_{\partial \Omega} \phi_- \cdot \partial_{\nu} \phi_-~.
\]
If both integrals in (\ref{eqn:boundary-integrals}) were negative, 
then the minimax characterization of the second eigenvalue would be 
contradicted.
\end{proof}

Let $e$ be a side of $T$.
Let $\partial_{\nu}$ be the unit outward normal derivative along $e$, and 
let $\partial_{\tau}$ denote the unit vector field tangent to 
$\partial P$ so that the frame $(\partial_{\nu}, \partial_{\tau})$ 
has the same orientation as the standard frame $(\partial_x, \partial_y)$.
In particular, $\partial_{\tau}$ corresponds to the counter-clockwise
orientation of $\partial P$. 

As far as we know, the following was first 
observed by Terence Tao in a post to \cite{Polymath}.

\begin{lem} \label{lem:boundary-integral-formula}
Let $L$ be a constant vector field. Suppose that $p$ (resp. $q$)
is a point on $e$ which is either a critical point or 
a vertex. If $\gamma:[a,b] \to e$, $\gamma(a)=p$, $\gamma(b)=q$,
and $\gamma' = \partial_{\tau}$,
then   
\[ 
\int_{\gamma}  Lu  \cdot \partial_{\nu} L u~
=~ 
-\frac{\mu}{2} 
\cdot 
\langle L, \partial_{\tau} \rangle 
\cdot 
\langle L, \partial_{\nu} \rangle
\cdot 
\left(
u(q)^2 - u(p)^2 
\right)~
\]
\end{lem}

\begin{proof}
Abusing notation slightly, we let $\partial_{\nu}$ (resp. $\partial_{\tau}$) 
denote the constant vector field that coincides with 
the restriction of $\partial_{\nu}$ (resp. $\partial_{\tau}$) to $e$.

We have 
$ 
L~ 
=~ 
\langle L, \partial_{\tau} \rangle \cdot \partial_{\tau}~
+~
\langle L, \partial_{\nu} \rangle \cdot \partial_{\nu}.
$
Thus, because $\partial_{\nu}u =0$ along $e$ and 
$\partial_{\nu} \partial_{\tau} = \partial_{\tau} \partial_{\nu}$,
we find that along $e$ we have 
$Lu= \langle L, \partial_{\tau} \rangle \cdot \partial_{\tau} u$
and 
$\partial_{\nu} Lu = \langle L, \partial_{\nu} \rangle \cdot \partial_{\nu}^2u$.
Hence
\[
Lu \cdot \partial_{\nu} Lu~
=~ 
\langle L, \partial_{\tau} \rangle 
\cdot 
\langle L, \partial_{\nu} \rangle 
\cdot 
\partial_{\tau} u
\cdot \partial_{\nu}^2 u.
\]
Since $\Delta u = \mu \cdot u$ we have 
$-\partial_{\nu}^2 u = \partial_{\tau}^2 u + \mu \cdot u$, and so 
\[
\frac{1}{C}
\cdot 
Lu \cdot \partial_{\nu} Lu~
=~  
-\partial_{\tau} u
\cdot 
\left(
\partial_{\tau}^2 u~
+~
\mu \cdot u
\right).
\]
where $C= \langle L, \partial_{\tau} \rangle 
\cdot 
\langle L, \partial_{\nu} \rangle$.
Since $\gamma'= \partial_{\tau}$
\begin{eqnarray*}
\frac{-1}{C}
\int_{\gamma}  
Lu \cdot \partial_{\nu} Lu~
&=&
\frac{1}{2}
\int_{\gamma} \partial_{\tau} \left((\partial_{\tau} u)^2 \right)~
+~ 
\frac{\mu}{2}
\int_{\gamma}  \partial_{\tau}(u^2) \\
&=&
\frac{1}{2}
\left( (\partial_{\tau} u(q))^2~ 
-~
(\partial_{\tau} u(p))^2\right)
+~ 
\frac{\mu}{2}
\left( u(q)^2~ 
-~
u(p)^2 \right).
\end{eqnarray*}
The term involving $\partial_\tau$ in the above expression vanishes 
if $p$ and $q$ are critical points of $u$.
If $p$ is a vertex then, using polar coordinates $(r, \theta)$ and the Bessel expansion of $u$ about $p$, 
a straightforward computation gives
\[
\partial_\tau u \left(r e^{i \th} \right)~
=~
a \cdot r \cdot \sin(\psi-\theta)~ +
o(r),
\]
if either the angle at $p$ is acute and $u(p) \neq 0$ or the angle at $p$ is obtuse and $c_1(p) = 0$, and 
\[
\partial_\tau u \left(r e^{i \th} \right)~
=~
a \cdot \nu \cdot r^{\nu-1} \cdot \cos(\psi - \th +\nu \cdot \th)~
+~
o(r^{\nu-1}),
\]
if either the angle at $p$ is obtuse and $c_1(p) \neq 0$ or the angle at $p$ is acute and $u(p) =0$.
If $q$ is a vertex as well then we have similar expressions about $q$ also.
The claim follows.
\end{proof}


\section{Proof of Theorem \ref{thm:main}}

The proof of Theorem \ref{thm:main} is by contradiction. We suppose that there exists 
an acute triangle $T_0$ such that a second Neumann eigenfunction 
$u_0$ of $T_0$ has an interior critical point. 
Let $T_t$, $0\leq t \leq 1$, be a continuous 
path of triangles that joins $T_0$ to a right isosceles triangle $T_1$
so that $T_t$ is acute for each $t<1$.
By the discussion in \S 12 of \cite{Annals}, for each $s \in [0,1]$ there exists 
a second Neumann eigenfunction $u_s$ so that $s \mapsto  u_s$ is continuous.

Recall that we call a critical point $p$ of $u_s$ {\it stable} if and only 
if for each neighborhood $U$ of $p$, there exists $\epsilon > 0$ 
so that if $|s - t| < \epsilon$ then $U$ contains a critical point of $u_t$.

\begin{lem}\label{lem:unstable-indexzero}
If $p$ is a critical point of $u_s$ with nonzero index,
then for each neighborhood $U$ of $p$, there exists $\epsilon>0$ 
such that if $|t-s|< \epsilon$ then $u_t$ has a critical point in $U$
whose index is the same as the index of $p$.
In particular, $p$ is a stable critical point.
\end{lem}

\begin{proof}
Let $\tp$ be a point in $DT$ that corresponds to $p$. Let $G_{\tp}$ be the 
level set of the lift $\tu_s$ that contains $\tp$.
Let $U$ be a disc neigborhood of $\tp$ that contains no other critical
points of $\tu_s$, and let $\gamma: S^1 \to \partial U$ be a parameterization 
of its boundary. The index of $\tp$ equals the winding number of 
$\nabla \tu_s \circ \gamma$ about $0 \in \Rbb^2$. Since $t \to \nabla u_t$ depends 
continuously on $t$, there exists $\epsilon>0$ such that if $|t-s|< \epsilon$,
then $\nabla \tu_t$ does not vanish on $\partial U$. 
In particular, the winding number of $\nabla \tu_t \circ \gamma$ about 0
is a well-defined and continuous function of $t \in (s-\epsilon, s+ \epsilon)$. 
Since this function is integer valued, it is constant and thus nonzero
by hypothesis. On the other hand, if $\tu_t$ had no critical point lying in $U$, 
then the winding number of $\nabla \tu_t \circ \gamma$ would equal zero. 
By Proposition \ref{prop:index-at-least-minus-one}, the 
index of each critical point of $\tu_s$ (resp. $\tu_t$) is either $+1$, $0$, or $-1$.
Thus, because the index is additive, it follows that
there is a critical point of $\tu_t$ in $U$ with the same index as $p$.
\end{proof}

A second Neumann eigenfunction of a right isosceles triangle has no 
critical points (see formula (26) in \cite{Annals}). Define 
$$
\tau~ 
:=~ 
\inf 
\left\{ 
t \in [0, 1]:\, u_t \mbox{ has at most two nonzero index critical points}  
\right\}.
$$
In particular, by Proposition \ref{prop:at-least-four},
we have $0< \tau \leq 1$. 

\begin{prop} \label{prop:one-or-two-stable}
The eigenfunction $u_{\tau}$ has no interior critical points, 
and $u_{\tau}$ has either one or two nonzero index critical points.
\end{prop}

\begin{proof}
The definition of $\tau$ and 
Lemma \ref{lem:unstable-indexzero} together imply that
the eigenfunction $u_{\tau}$ has at most two nonzero 
index critical points. In particular, by Proposition 
\ref{prop:at-least-four}, 
the function $u_\tau$ has no critical points in the interior of $T_\tau$.

By the definition of $\tau$, there exists a sequence $t_n \nearrow \tau$ 
such that $u_{t_n}$ has at least three nonzero index critical points.  
Without loss of generality, by Proposition 7.10 \cite{Annals}, we may assume 
that for each $n$, two of these critical points, $p_n$ and $q_n$,
lie on different sides of $T_{t_n}$. 

We claim that one of the sequences has a limit point which is not a vertex. 
Indeed, first note that by Lemma 9.2 \cite{Annals}, the sequences $p_n$ 
and $q_n$ do not 
have a common limit point which is a vertex. Next, suppose that a vertex $v$ 
is a limit point of $p_n$ and a vertex $v' \neq v$ is a limit point of $q_n$.
If $\tau<1$, then the function $u_{\tau}$ 
vanishes at two vertices thus contradicting Corollary 5.5 \cite{Annals}.
If $\tau=1$, then $T_{\tau}$ is an isosceles right triangle 
and so either $v$ or $v'$ is an acute vertex. Thus, by Corollary 9.3 \cite{Annals},
the function $u_{\tau}$ vanishes at this vertex. But inspection of 
the explicit form of $u_{\tau}$ on a right isosceles triangle shows
that $u_{\tau}$ does not vanish at an acute vertex. 

Let $z$ be the limit of one of the two sequences that lies in
the interior of a side $e$. Then $z$ is a critical point of $u_\tau$.
If the index of the critical point $z$ of $u_{\tau}$ equals zero, then Corollary \ref{coro:n-arcs-n-1-stable} implies that $u_{\tau}$ has a 
nonzero index critical point in $\partial T \setminus e$.
\end{proof}

\begin{figure}
\includegraphics[scale=.5]{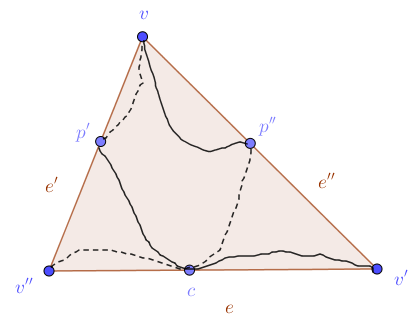} 
\caption{The nodal set of $L_{e'} u_{\tau}$ (solid) 
and the nodal set of $L_{e''} u_{\tau}$ (dashed) in the proof of 
Proposition \ref{prop:exactly-one-stable}. \label{fig:two-stable-two}}
\end{figure}

\begin{prop} \label{prop:exactly-one-stable}
The eigenfunction $u_{\tau}$ has exactly
one nonzero index critical point, and this critical point
has index equal to $-1$. 
Moreover, $u_\tau$ does not vanish at any vertex of 
$T_\tau$.
\end{prop}

\begin{proof}
By Proposition \ref{prop:one-or-two-stable}, the function $u_\tau$ has 
either one or two nonzero index critical points.  By Proposition
\ref{prop:index-at-least-minus-one}, each such critical point has index 
equal to $+1$ or $-1$. 

By Lemma \ref{lem:pos-local-max}, a vertex $v$ is a local extremum 
if and only if $u_\tau(v) \neq 0$. By Corollary 5.5 \cite{Annals}, 
the function $u$
vanishes at at most one vertex, and hence either two or three vertices 
of $T_\tau$ are local extrema. Each such local extremum has index
equal to $+1$.  Thus, it follows from Proposition \ref{prop:index-thm}
that if $u_\tau$ vanishes at a vertex, then $u_\tau$ has two nonzero index
critical points, and if $u_\tau$ does not vanish at a vertex, then 
$u_\tau$ has exactly one nonzero index critical point and this critical
point has index $-1$.  

Thus, to prove the proposition, we need only show that $u_{\tau}$ 
does not vanish at a vertex. Suppose to the contrary that $u_\tau(v)=0$ 
for some vertex $v$ and hence $u_\tau$ has exactly two nonzero index 
critical points. 
In particular, the nodal set $\Zcal(u_\tau)$ is a simple arc with one endpoint $v$. 
By Lemma 3.3 \cite{Annals}, the other endpoint of $\Zcal(u_\tau)$ lies 
in the interior of  the side $e$ opposite to $v$.  
Let $v'$, $v''$ denote the other two vertices of $T$, and 
let $e'$ and $e''$ denote, respectively, the opposite sides.


Note that no side contains both of the nonzero index critical points. 
Indeed, if, for example, the side $e$ were to contain both 
of these critical points then,
by Corollary  \ref{coro:n-arcs-n-1-stable}, 
we would have at least
two arcs of $\Zcal(L_e u_{\tau})$ that emanated from $e$ into the interior
of $T$. Lemma \ref{prop:degree-1-stable} would then imply the existence 
of a third  critical point, contradicting 
Proposition \ref{prop:one-or-two-stable}. 

The values of $u_\tau$ have opposite signs on the two components of 
$T_{\tau} \setminus \Zcal(u_\tau)$. Hence,  
because $v'$ and $v''$ lie in distinct components,
we have $u_{\tau}(v') \cdot u_{\tau}(v'') < 0$. Thus, by Proposition 
\ref{prop:crit_opposite_sign_on_vertices}, the number of critical points 
with non-zero index that lie on $e$ is even and hence equals zero 
by the previous paragraph.
In sum, one nonzero index critical point, $p'$, lies on $e'$ and the other, 
$p''$, lies on $e''$.

Since $u(v)=0$ and $p'$ is a nonzero index critical point,
Corollary  \ref{coro:n-arcs-n-1-stable} implies that  
$\Zcal(L_{e'} u_\tau) \setminus e'$ has at least two components each of 
whose closures intersect $e'$, each of these
component has a degree 1 vertex in $\partial T \setminus e'$, 
and the closure of one of these components
contains the vertex $v$ and the closure of the other component
contains $p'$.
Because $p'$ and $p''$  are the only nonzero index
critical points, by Proposition \ref{prop:degree-1-stable}, 
the points $v'$ and $p''$ are the only
possible degree 1 vertices of 
$\mathcal{Z}(L_{e'} u_\tau) \setminus e'$ 
that lie in $\partial T \setminus e'$. Hence one 
component of $\Zcal(L_{e'} u_\tau) \setminus e'$ is an arc
$\alpha'$ that joins $p'$ to $v'$ and the other component is an
arc that joins $v$ to $p''$.

A symmetric argument shows that an arc component $\alpha''$ of 
$\Zcal(L_{e''} u) \setminus e''$ joins $p''$ to $v''$ and the other 
component is an arc that joins $v$ to $p'$. Since $\alpha'$ 
separates $p''$ from $v''$, the arcs $\alpha'$ and $\alpha''$ 
intersect at some point $c$.
If $c$ were to lie in the interior of $T$, then it would be 
an interior critical point contradicting the fact that $u_{\tau}$
has no interior critical points. 

Suppose that $c$ were to lie on the boundary of $T$.
Since $p'$ and $p''$ are nonzero index critical points, the intersection
point $q$ does not equal $p'$ or $p''$. 
Since $\alpha'$ and $\alpha''$ are arc components of 
$\Zcal(L_{e'} u_\tau) \setminus e'$ and $\Zcal(L_{e''} u_\tau) \setminus e''$
respectively, the intersection point $c$ does not lie on $e'$ or $e''$. 
Hence $c$ lies in the interior of the side $e$, and hence is a
critical point of $u_\tau$, and, moreover, 
$c$ is a degree two vertex of both $\Zcal(L_{e'}u)$ 
and $\Zcal(L_{e''}u)$. See Figure \ref{fig:two-stable-two}.
Since $T$ is acute, this contradicts Proposition \ref{prop:degree-two-implies-degree-zero}.
%
%
%
\end{proof}

\begin{figure}
\includegraphics[scale=.5]{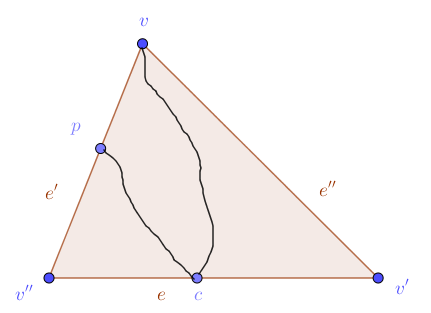}
\caption{The nodal set $L_e u$ in Lemma \ref{lem:cusp-two-components}. \label{fig:one-stable-cusp}}
\end{figure}

\begin{lem} \label{lem:cusp-two-components}
Let $e'$ be the side that contains the critical point 
$p$ of $u_\tau$ 
with index equal to $-1$. Suppose that $c$ is an index zero critical point of $u_{\tau}$ that lies on a side $e \neq e'$. 
Then $\Zcal(L_e u_\tau) \setminus e$ has exactly  
two components. One component is an arc that joins $c$ to $p$, 
and the other is an arc that joins $c$ to the vertex $v$ opposite to $e$.
\end{lem}

\begin{proof}
By Corollary  \ref{coro:n-arcs-n-1-stable}, the set $\Zcal(L_e u_\tau) \setminus e$ 
has at least two components each of whose closures 
$C$ and $C'$ intersect $e$, and, moreover, $C$ (resp. $C'$) has a 
degree 1 vertex that lies in $\partial T \setminus e$ and the only degree 1 vertex
of $C$ (resp. $C'$) that lies in $e$ is $c$. By Proposition \ref{prop:degree-1-stable},
the only possible degree 1 vertices are the nonzero index critical
point $p$ and the vertex $v$ opposite to $e$.
Hence, by relabeling if necessary, we may assume
that $C$
is an arc that joins $c$ to $p$ and $C'$ is an arc that 
joins $c$ to $v'$. 

If $\Zcal(L_e u_\tau) \setminus e$ 
were to have another component, then the eigenfunction 
$u_\tau$ would have a nonzero index critical point in $\partial T \setminus e$
that would be distinct from $p$. But $p$ is the only nonzero index critical point.
\end{proof}

\begin{prop} \label{prop:saddle-and-cusp}
The eigenfunction $u_\tau$ has exactly two critical points and
these critical points lie on distinct sides of $T_\tau$.  
One critical point has index equal to $-1$, and the other critical 
point has index equal to zero.
\end{prop}

\begin{proof}
By Proposition \ref{prop:exactly-one-stable}, the function $u_\tau$ 
has exactly one nonzero index critical point $p$, the critical point 
$p$ is a saddle, and, moreover, $u_\tau$ does not vanish at a vertex.
In particular, it suffices to 
show that $u_\tau$ has exactly one  critical point of index zero.
Let $e'$ be the side of $T_{\tau}$ that contains the unique 
nonzero index critical point $p$. 

First we note that $p$ is the only 
critical point on $e'$. Indeed, suppose that a second critical point $q$ 
belonged to $e'$. Then, by Corollary \ref{coro:n-arcs-n-1-stable}, the
set $\Zcal(L_{e'} u) \ e$ would have at least three degree 1 vertices
in $\partial T \setminus e'$. Except for the vertex $v'$ opposite to $e'$,
each degree 1 vertex in $\partial T \setminus e'$ 
is a nonzero index critical point 
by Proposition \ref{prop:degree-1-stable}. Thus we would
have two nonzero index critical points, thus contradicting 
Proposition \ref{prop:exactly-one-stable}.

A similar argument shows that $u_{\tau}$ cannot
have more than one critical point of index zero on some side of $T_\tau$. 
Thus, to finish the proof, it suffices to rule out the possibility that $u_\tau$ 
has both an  index zero critical point $c$ on $e$ and 
an  index zero critical point $c''$ on $e''$.

By Lemma \ref{lem:cusp-two-components}, the set 
$\mathcal{Z}(L_e u_\tau) \setminus e$
contains an arc $\alpha$ that joins $c$ and $p$ and an arc $\beta$ 
that joins $c$ to $v$. Similarly, the set 
$\mathcal{Z}(L_{e''} u_\tau) \setminus e''$
consists of an arc $\alpha''$ that joins $c''$ and $p$, and an arc 
$\beta''$ that joins $c''$ to $v''$.

If $\beta$ (resp. $\beta''$) does not intersect the interior of $e''$ at $c''$
(resp. $e$ at $c$), 
then $\beta$ (resp. $\beta''$) intersect $\alpha''$ (resp. $\alpha$)
in the interior of $T$, and hence we have an interior critical point, 
a contradiction.
Thus, we may assume that $\beta$ (resp. $\beta''$) intersects the interior of 
$e''$ at $c''$ (resp. $e$ at $c$). Since $c$ (resp. $c''$) is the
only critical point of $u_\tau$ on $e$ (resp. $e''$)
the arc $\alpha$ (resp. $\alpha''$) has a degree 2 vertex at the 
 index zero critical point $c$ (resp. $c''$). See Figure \ref{fig:cusp_on_each_side}. 
 In particular, $c$ (resp. $c''$) is a degree two vertex of $\Zcal(L_e u_\tau)$ and $\Zcal(L_{e''} u_\tau)$ (resp. $\Zcal(L_{e''} u_\tau)$ and $\Zcal(L_e u_\tau)$).

\begin{figure}
\includegraphics[scale=.5]{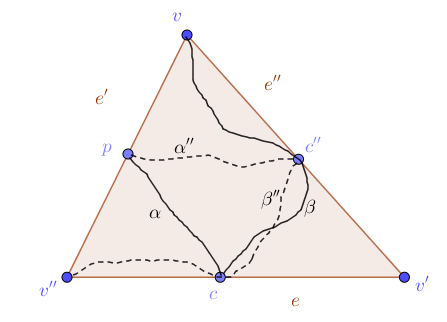}
\caption{The nodal sets of $L_e u_{\tau}$ (solid) and $L_{e''} u_{\tau}$ (dashes)
in the proof of Proposition \ref{prop:saddle-and-cusp}.
\label{fig:cusp_on_each_side}}
\end{figure}

By Proposition \ref{prop:degree-1-stable},
the set $\Zcal(L_{e'} u_\tau)$ has at least one component 
whose closure intersects $e'$ and this component has a degree 1 vertex 
that belongs to  $\partial T \setminus e'$. Since $p$ is
the only nonzero index critical point of $u_\tau$, the vertex $v'$ is the 
only possible degree 1 vertex of $\Zcal(L_{e'} u_\tau)$ 
that lies in $\partial T \setminus e'$, and hence 
$\Zcal(L_{e'} u_\tau)$ consists of an arc $\gamma$ in 
$\Zcal(L_{e'} u_\tau)$ that joins $p$ to $v'$. In particular, since $p$ 
and $v'$ lie in distinct components of $T \setminus \beta$, the arc $\gamma$
intersects $\beta$. Since $u_\tau$ has no interior critical points, 
the intersection does not occur in the interior of $T$, and so either 
$\gamma$ contains $c$ or $\gamma$ contains $c''$. By relabeling if necessary,
we may assume that $\gamma$ contains $c$.

In sum, the critical point $c$ is a degree 2 vertex of 
$\Zcal(L_e u_\tau)$ and $\Zcal(L_{e''} u_\tau)$.
Proposition \ref{prop:degree-two-implies-degree-zero} implies 
that $c$ is a degree 0 vertex of $\Zcal(L_{e'} u_\tau)$,
a contradiction.
\end{proof}

In what follows we suppose that the saddle critical point of $u_{\tau}$
lies on the side $e'$ and that the  index zero critical point lies on $e$.  

\begin{prop} \label{prop:u_t-description}
There exists $\epsilon>0$ so that if $t \in (\tau - \epsilon, \tau)$, then
the critical points of $u_t$ consist of
\begin{enumerate}

\item exactly two nonzero index critical points on $e$ and no other 
critical points on $e$, \item a saddle critical point on the side $e'$ 
and no other critical points on $e'$,

\item no critical points on  the side $e''$, and

\item no interior critical points.
\end{enumerate}
\end{prop}

\begin{proof}
Let $p$ be the saddle critical point of $u_\tau$ and let $q$ be the 
index zero critical point of $u_\tau$ given by Proposition
\ref{prop:saddle-and-cusp}.
Let $e'$ be the side that contains $p$, and let $e \neq e'$ be the side
that contains $q$. 

Because $p$ is a saddle critical point of $u_{\tau}$, by Lemma \ref{lem:unstable-indexzero}, there exists $\epsilon_1>0$ 
and a neighborhood $U$ of $p$ such that if $|t- \tau|<\epsilon_1$, then 
$u_t$ has exactly one critical point in $U$ and this critical point is a saddle. 
Because $u_{\tau}$ has no critical points on $e''$ and $u_{\tau}$ does not
vanish at a vertex of $T$, there exists $\epsilon_2>0$ such that if 
$|t-\tau| < \epsilon_2$, then $u_{\tau}$ has no critical point on $e''$ 
and hence no interior critical point by 
Proposition \ref{prop:at-least-four}.

We claim that there exists $\epsilon_3>0$ so that if $\tau - \epsilon_3 < t < \tau$,
then each critical point of $u_t$ other than $p$ belongs to
the side $e$. Indeed, otherwise there would be a sequence $t_n \nearrow \tau$ 
and critical points $q_n$ of $u_{t_n}$ that belong to either $e \setminus U$  
or to the side $e''$. By Proposition \ref{prop:one-or-two-stable}, the
function $u_{\tau}$ does not vanish at a vertex of $T$, and so, by Lemma 
9.3 \cite{Annals}, no vertex is an accumulation point of $q_n$.
Hence each accumulation point lies in the interior of a side.
But such an accumulation point would be a critical point of $u_\tau$,
and $p$ and $q$ are the only critical points of $u_{\tau}$.

By the definition of $\tau$, for each $t < \tau$, the function 
$u_t$ has at least three nonzero index  critical points. 
Since $u_t$ has no critical points on $e''$ and only one 
critical point on $e'$ for $t \in (\tau -\epsilon_3,\tau)$,
the function $u_t$ has at least two nonzero index critical points on $e$. 

To finish the proof it suffices to show that, 
for $t$ sufficiently near $\tau$,
the function $u_t$ has at most two critical points on $e$.
If there were more than two critical points on $e$, 
then there would exist a sequence $t_n \nearrow \tau$ and 
critical points $p_n$, $q_n$, and $r_n$ of $u_{t_n}$ that lie on $e$.
Since $u_\tau$ does not vanish at a vertex and has only one critical
point on $e$, each sequence would limit 
to the  index zero critical point $q$ of $u_\tau$.
For each $n$, Corollary  \ref{coro:n-arcs-n-1-stable} 
implies that the set 
$\Zcal(L_e u_{t_n}) \setminus e$ has at least 
three components whose closures intersect $e$.
By considering the limit we find that $\Zcal(L_e u_{\tau}) \setminus e$
would have at least three components whose closures intersect $e$. 
But the only possible degree 1 vertices of $\Zcal(L_e u_{\tau}) \setminus e$
that lie in $\partial T \setminus e$ are $p$ and $v$, a contradiction.
\end{proof}

Let $\epsilon$ be as in Proposition \ref{prop:u_t-description}, and 
define $\xi$ to be the supremum of $t$ such that $u_t$ has more than three
critical points, that is
$$ \xi~ :=~ \sup \left\{ t < \tau -\epsilon/2\, :\, 
u_t \mbox{ has at least 4 critical points} \right\}.$$
The assumption that $u_0$ has an interior critical 
point combined with Proposition \ref{prop:at-least-four} 
implies that $\xi \geq 0$.

\begin{prop} \label{prop:u-xi}
The function $u_{\xi}$ satisfies the following conditions:
\begin{enumerate}

\item $u_{\xi}$ has exactly two nonzero index critical points on $e$, 
\item $u_{\xi}$ has exactly one saddle critical point on the side $e'$, 
\item $u_{\xi}$ has no nonzero index critical points on the side $e''$, 
\item $u_{\xi}$ has no interior critical points,
\item $u_\xi$ has at least one  index zero critical point,  and
\item $u_{\xi}$ does not vanish at a vertex of $T_{\xi}$.
\end{enumerate}
\end{prop}

\begin{proof}
By definition of $\xi$ and Proposition \ref{prop:u_t-description}, 
we have $\xi < \tau$.
We claim that for $t \in [\xi, \tau)$, the eigenfunction $u_t$ has exactly
three nonzero index critical points. Indeed, since $t< \tau$, the function
$u_t$ has at least three nonzero index critical points. 
Hence, if $\xi < t < \tau$, then by the definition of $\xi$, the function $u_t$
has exactly three critical points and each has nonzero index. 
Thus, Proposition \ref{lem:unstable-indexzero} and the definition of $\tau$
imply that $u_{\xi}$ has exactly three nonzero index critical points.

Let $t \in [\xi, \tau)$ and let $z_1(t)$, $z_2(t)$, and $z_3(t)$
be the three nonzero index critical points of $u_t$. 
By Lemma \ref{lem:unstable-indexzero}, 
for each open set $U_i$ containing $z_i(t)$, there exists $\epsilon_i>0$
such that if $|s-t|  < \epsilon_i$, then there exists 
a nonzero index critical point $z_i(s)$ that lies in $U_i$. 
Suppose that the sets $U_i$ are mutually disjoint.
Then since $u_s$ has exactly three nonzero index critical points, the 
critical point $z_i(s)$ is the only nonzero index critical point
of $u_s$ that lies in $U_i$. Using Lemma \ref{lem:unstable-indexzero} 
we find that $s \mapsto z_i(s)$ is continuous on 
$(t-\epsilon^*, t+ \epsilon^*)$
where $\epsilon^* = \min \{\epsilon_1, \epsilon_2,\epsilon_3\}$.
One may cover $[\xi, \tau - \epsilon/2]$ with 
finitely many intervals of this form.
The fact that $u_s$ has exactly three nonzero index critical points 
and Proposition \ref{prop:u_t-description}
imply that we may parameterize these critical points
with continuous functions
$s \mapsto z_i(s)$ defined on all of $[\xi, \tau)$.

Thus, since the index is continuous by Lemma \ref{lem:unstable-indexzero},
conditions (1), (2), and (3) hold. 
Condition (4) follows from Proposition \ref{prop:at-least-four} 
and condition (3).  
Condition (5) follows from the definition of $\xi$. 

By Corollary 5.5 \cite{Annals}, the eigenfunction $u_\xi$ can vanish at at most 
one vertex. By Proposition \ref{prop:index-at-least-minus-one},
the index of each of the three nonzero index critical points of $u_{\xi}$
equals $\pm 1$. Hence  condition (6) follows from Proposition \ref{prop:index-thm}.
\end{proof}

To prove Theorem \ref{thm:main} we need only prove that, 
for each acute triangle $T$, a second Neumann eigenfunction $u$ of $T$ 
that satisfies conditions (1) through (6) in Proposition \ref{prop:u-xi} cannot exist. 
We suppose to the contrary the existence of such an eigenfunction $u$ and then 
deduce a contradiction.

Let $p$ denote the unique nonzero index critical point on $e'$.
Let $q$ and $q'$ denote the nonzero index critical points on $e$ 
labeled so that $q$ lies on the segment $[v''q']$.

We will proceed in three steps. 
In step 1, we will show that $u$ does not have 
an index zero critical point $c$ on $e$, and we will use
this to show that the critical points $q$ and $q'$ are nondegenerate.
In step 2, we will show that $u$ does not have 
an index zero critical point $c$ on $e'$, and we will use
this to show that the critical point $p$ is nondegenerate.
Finally, in step 3, we will use the nondegenracy of $p$, $q$, and $q'$
to conclude that $u$ has no index zero critical points on $e''$.
Hence $u$ has no zero index critical points, thus contradicting condition (6).

\hspace{.1cm}
\paragraph{\bf Step 1: $c$ does not belong to $e$ and $q$, $q'$ are nondegenerate.}

Suppose to the contrary that
the index zero critical point $c$ belongs to the side $e$. Then,
by Corollary  \ref{coro:n-arcs-n-1-stable}, and 
Proposition \ref{prop:degree-1-stable}, the set  
$\mathcal{Z}(L_e u) \setminus e$ has at least three nonzero index   
critical points in $\partial T \setminus e$
contradicting conditions (2) and (3) in Proposition \ref{prop:u-xi}.
In particular, the only critical points of $u$ that lie on $e$
are the critical points $q$ and $q'$ by condition (1).

Suppose that $q$ were degenerate.   
By Lemma 7.2 \cite{Annals}, the eigenfunction $u$ has the form
\[ 
u(z)~ 
=~
a_{00}~ +~ a_{20} \cdot x^2~ +~ a_{02} \cdot y^2~ +~ 
3 a_{30} \cdot (x^3 -3 x\cdot y^2)~
+~
o(|z|^3)
\]
near $q$.
Since $q$ is degenerate either $a_{20}$ or $a_{02}$ equals zero.
If $a_{02}=0$ then by  Proposition 7.3 \cite{Annals} 
and Proposition \ref{prop:degree-1-stable}
there would be a nonzero index critical point in each side of $T$.
But this would contradict condition (3) in Proposition \ref{prop:u-xi}. 
Hence $a_{20}=0$ and 
$\partial_x u(z) =  3 a_{30} \cdot ( x^2 -y^2) + o(|z|^2)$. 
In particular, there are at least two arcs of $\Zcal(L_e u)$ 
emanating from $q$. By Lemma 6.6 \cite{Annals}, at least one arc of 
$\Zcal(L_e u)$ emanates from $q'$. Thus,
the set $\mathcal{Z}(L_e u) \setminus e$ 
would have at least three  components whose closures intersect $e$ 
and hence $\mathcal{Z}(L_e u) \setminus e$ would have at 
least two nonzero index critical points in $\partial T \setminus e$
contradicting conditions (2) and (3) in Proposition \ref{prop:u-xi}.
Thus, $q$ can not be degenerate, and by the same argument $q'$ is
also nondegenerate.

\hspace{.1cm}
\paragraph{\bf Step 2: $c$ does not belong to $e'$ and $p$ is nondegenerate.} 

Suppose to the contrary that the  index zero 
critical point $c$ belongs to the side $e'$.

Since $e$ contains two nonzero index critical points, 
Corollary   \ref{coro:n-arcs-n-1-stable} implies
the set  $\mathcal{Z}(L_{e} u) \setminus e$ has  
at least two components whose closures intersect $e$, and each component has 
a degree one vertex in $\partial T \setminus e$.
By Proposition \ref{prop:degree-1-stable} and Proposition \ref{prop:u-xi}, the only possible degree 1 vertices in $\partial T \setminus e$ 
are $v$ and $p$, and so one component is an arc that joins 
$q$ to $p$ and the other is an arc $\delta$ that joins $q'$ to $v$.

Since $c$ is a critical point of index zero, 
the proof of Corollary  \ref{coro:n-arcs-n-1-stable} implies
that at least two arcs of $\mathcal{Z}(L_{e'} u) \setminus e'$ 
emanate from $c$ and at least one arc emanates from $p$. 
The only possible degree 1 vertices of $\mathcal{Z}(L_{e'} u) \setminus e'$ 
in $\partial T \setminus e'$ are the vertex $v'$ and 
the nonzero index critical points $q$ and $q'$ on $e$. Hence there 
are exactly three arc components and each
arc joins either $p$ or $c$ to either $q$, $q'$, or $v$.

By arguing as in step 1 above, we can further conclude 
that $p$ is non-degenerate and $e'$ contains at most one zero index critical point.
 
The  index zero critical point 
$c$ belongs to either the line segment $[v''p]$ or the segment $[pv]$.
We claim that $c$ does not belong to the line segment $[v''p]$.
Indeed, suppose to the contrary that $c$ belongs to the line segment $[v''p]$.
Then there exists an arc component $\beta$ of $\Zcal(L_{e'} u) \setminus e'$
that joins $p$ to $v'$. 
The other components of $\Zcal(L_{e'} u) \setminus e'$ are arcs that 
join $c$ to $q$ and $q'$ respectively.
 Note that the arc $\beta$  necessarily 
intersects the arc $\delta$ that joins $q'$ to $v$.
Since $u$ has no interior critical points, the intersection does not
occur in the interior of $T$.  The arc $\beta$ (resp. $\delta$) is 
disjoint from the component of $\Zcal(L_e u)$  (resp. $\Zcal(L_{e'} u)$) 
that contains $p$, $c$ and $q$  (resp. $c$, $q$, and $q'$).
It follows that the intersection point does not occur in $\partial T$.
Thus, we have a contradiction, and $c$ does not lie on the line segment $[v''p]$.

\begin{figure}
\includegraphics[scale=.5]{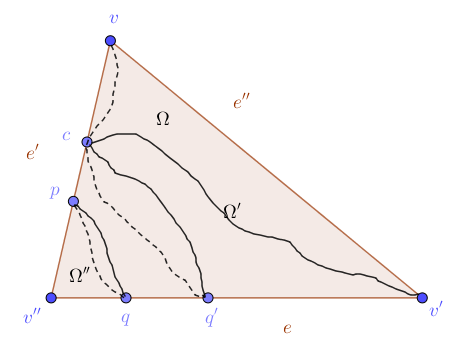}
\caption{Step 2. The nodal sets $\Zcal(L_e u)$ (dashed) and  $\Zcal(L_{e'} u)$
(solid).
\label{fig:second-case_c_on_right.png}}

\end{figure}

Thus, we may assume that $c$ lies on the line segment $[pv]$.
In this case, one component of  $\mathcal{Z}(L_{e'} u) \setminus e'$ joins 
$p$ to $q$, a second component joins $c$ to $q'$, and the third joins $c$ to $v'$. 
Since the arc $\delta$ that joins $q'$ to $v$ does not intersect 
$\Zcal(L_{e'} u)$ in the interior of $T$, the arc $\delta$ contains $c$.
See Figure \ref{fig:second-case_c_on_right.png}. 

Let $\Omega$ be the component of  $T \setminus \Zcal(L_{e'} u)$ that contains $v$,
let $\Omega'$ be the component of  $T \setminus \Zcal(L_e u)$ that contains $v'$,
and let $\Omega''$ be the component of $T \setminus \Zcal(L_e u)$ that contains $v''$.
Note that $\Omega \cap \Omega'' = \emptyset$ and $\Omega' \cap \Omega'' = \emptyset$.
Thus, to achieve a contradiction in this case, by Proposition
\ref{lem:boundary-term-nonneg},
it suffices to show that the boundary integral 
of $L_e u$ over $\Omega''$ is negative, and that either the boundary 
integral of $L_e u$ over $\Omega$ is negative or the boundary integral 
of $L_{e'} u$ over $\Omega'$ is negative.

Using Lemma \ref{lem:boundary-integral-formula}, we find that 
the sign of the boundary integral of $L_e u$ over $\Omega''$ is 
opposite to the sign of $u(v'')^2-u(p)^2$. Hence to show that this
boundary integral is negative it suffices to show that $|u(v'')| > |u(p)|$. 
Since $e'$ has exactly one nonzero index critical point and 
$u$ does not vanish at the endpoints of $e'$,
Proposition \ref{prop:crit_opposite_sign_on_vertices} implies that $u$
does not vanish on $e'$. 
By Lemma \ref{lem:pos-local-max}, the vertex $v$ is either a 
positive local maximum or a negative local minimum. Because no nonzero index 
critical point lies between $p$ and $v''$, using 
Lemma \ref{lem:tangential-cusp} we deduce that $|u(v'')| > |u(p)|$
as desired. Thus, the boundary integral of 
$L_e u$ over $\Omega''$ is negative.

Since $u$ satisfies Neumann boundary conditions, $e''$ is the only side of $T$ over which we 
need to compute the 
boundary integral of both $L_{e'}u$ over $\Omega$ and $L_{e}u$ over $\Omega'$. 
Since $T$ is acute, the products 
$\langle L_{e}, \partial_{\tau} \rangle \cdot \langle L_{e}, \partial_{\nu} \rangle$ 
and 
$\langle L_{e'}, \partial_{\tau} \rangle \cdot \langle L_{e'}, \partial_{\nu} \rangle$ 
have opposite signs. Thus, it follows from Proposition 
\ref{lem:boundary-integral-formula} that either the boundary integral of
$L_{e'}u$ over $\Omega$ is negative or the boundary integral of 
$L_{e}u$ over $\Omega'$ is negative.

\hspace{.1cm}
\paragraph{\bf Step 3: $c$ does not belong to $e''$.} 

Suppose to the contrary that the  index zero critical point $c$ belongs to the side $e''$.
As in  step 2, Corollary  \ref{coro:n-arcs-n-1-stable}  
implies that one component 
of $\mathcal{Z}(L_{e} u)$ is an arc joining $q$ to $p$ and the other 
component is an arc joining $q'$ to $v$. In this case, 
the arc joining $q$ to $p$ does not pass through $c$. 

Let $\Omega''$ be the component of $T \setminus \Zcal(L_e u)$ that contains $v''$. 
As observed in step 2, $u$ does not vanish on $e'$ and hence $|u(v'')|>|u(p)|$.
Thus it follows from Lemma \ref{lem:boundary-integral-formula} that 
the boundary integral of $L_e u$ over $\Omega''$ is negative.
Hence, by Lemma \ref{lem:boundary-term-nonneg},
to obtain a contradiction, it suffices to find a domain disjoint from $\Omega''$
over which a boundary integral is negative.

\begin{figure}
\includegraphics[scale=.5]{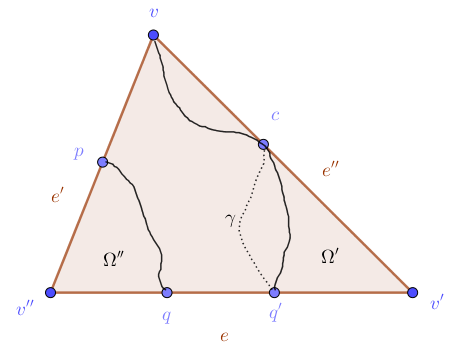}
\caption{Step 3: The nodal set of $L_e u$ (solid).
\label{fig:L_e_case_3}} 
\end{figure}

Consider the nodal set $\mathcal{Z}(L_{e''} u)$.
Since $c$ has index zero, at least two components of 
$\mathcal{Z}(L_{e''} u) \setminus e''$ have closure containing $c$.
By steps 1 and 2, the critical points $p$, $q$ and $q'$ 
are nondegenerate, and hence each is a degree 1 vertex of 
$\mathcal{Z}(L_{e''} u)$. Since $u(v'') \neq 0$, the vertex $v''$
is also a degree 1 vertex. Since $u$ has no interior critical points, the sets
$\mathcal{Z}(L_{e} u)$ and $\mathcal{Z}(L_{e''} u)$ do not intersect
in the interior of $T$. It follows that the arc of $\Zcal(L_e u)$
that joins $q'$ to $v$ contains the  index zero critical point $c$.
Also, it follows that the component of 
$\mathcal{Z}(L_{e''} u)$ that contains $v''$ is an arc whose 
other endpoint is $p$ or $q$. On the other hand, 
by Corollary  \ref{coro:n-arcs-n-1-stable}, 
the set $\mathcal{Z}(L_{e''} u) \setminus e''$ 
has at least two components whose closures intersect $e''$.
It follows that there exists an arc $\gamma$
contained in $\Zcal(L_{e''}u)$ that joins $c$ and $q'$. 
Figure \ref{fig:L_e_case_3}.

Because $u(v') \neq 0$, the vertex $v'$ is a degree 1 vertex of 
$\Zcal(L_{e'}u)$. Since $\Zcal(L_{e'}u)$ and $\Zcal(L_{e''}u)$ do not intersect 
in the interior of $T$, the component of $\Zcal(L_{e'}u)$ which contains $v'$
either is an arc $\alpha$ that joins $v'$ to $q'$ or it contains an 
arc $\alpha$ that joins $v'$ to $c$. 
We consider these two cases separately.

\begin{figure}
\includegraphics[scale=.5]{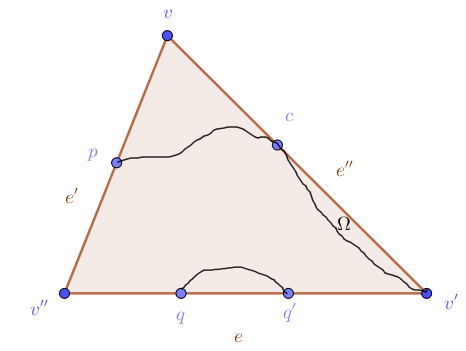}
\caption{Step 3: The nodal set of $L_{e'}u$  (solid) when $\alpha$
joins $v'$ to $c$. \label{fig:L_e_case_3a}}
\end{figure}

If $\alpha$ joins $v'$ to $c$, then let 
$\Omega$ be the component of $T \setminus \alpha$ that is 
bounded by $\alpha$  and the segment $[cv']$. See Figure \ref{fig:L_e_case_3a}. 

By Proposition \ref{lem:boundary-integral-positive}, the boundary integral 
of $L_{e'} u$ over $\Omega$ is positive. 
Hence it follows from Lemma \ref{lem:boundary-integral-formula}
that $u(c)^2 -u(v')^2$ is negative.
On the other hand, if $\Omega'$ is the component of $T \setminus \Zcal(L_e u)$ 
that contains $v'$, then by Lemma \ref{lem:boundary-integral-formula}, the 
boundary integral of $L_eu$ over $\Omega'$ has the same sign as 
$u(c)^2 -u(v')^2$, and so it is negative. 
Since $\Omega''$ and $\Omega'$ are disjoint, this contradicts 
Proposition \ref{lem:boundary-term-nonneg}.

If $\alpha$ joins $v'$ to $q$, then let 
$\Omega'''$ be the component of $T \setminus \Zcal(L_{e''} u)$
that contains $v'$. To exhibit a contradiction it suffices to show that 
the boundary integral of $L_{e''} u$ over $\Omega'''$ is negative. 
To see this, note that 
the boundary integral of $L_{e'}u$ over the component 
of $T \setminus \Zcal(L_{e'}u)$ that contains $[q'v']$ is
positive by Proposition \ref{lem:boundary-integral-positive}.
But Proposition \ref{lem:boundary-integral-formula}
implies that the sign of this boundary integral 
is opposite to the sign of the boundary  integral of $L_{e''}u$ 
over $\Omega'''$.

\begin{figure}
\includegraphics[scale=.5]{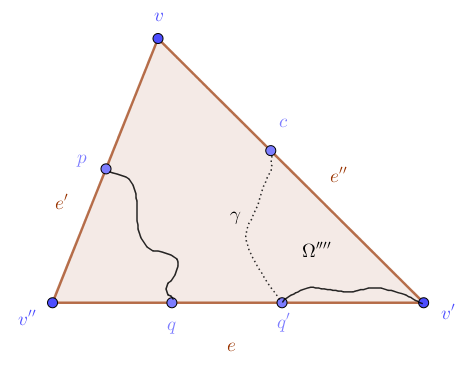}
\caption{Step 3: The nodal set of $L_{e'}u$ (solid)
when $\alpha$ joins $v'$ to $q$.
The arc $\gamma$ (dotted) lies in $\Zcal(L_{e''}u)$. 
\label{fig:L_e_case_3b}} 
\end{figure}


\section{The cases of obtuse and right triangles}
\label{sec:obtuse-right}

As indicated in the introduction, the proof 
of Theorem 1.1 in \cite{Annals} for both obtuse and right triangles
follows from the method presented there.  Namely, we have the following.

\begin{thm} \label{thm:obtuse-right}
If $T$ is an obtuse or right triangle and  
$u$ is a second Neumann eigenfunction on $T$, 
then $u$ has no critical 
points.\footnote{As in \cite{Annals}, we do not consider a local extremum 
at a vertex to be a possible critical point.}
\end{thm}

To prove Theorem \ref{thm:obtuse-right} in the obtuse case we use
the following modifications of Proposition 7.6 and Theorem 7.7 \cite{Annals}.
Note that Propositions \ref{prop:tangential-cusp} and \ref{prep:thm} below
hold true for a second Neumann eigenfunction $u$ on any triangle. 
Lastly, we would like to point out that, although the statement of Theorem \ref{thm:obtuse-right} is stronger than Theorem 1.1 in \cite{Annals}, Theorem \ref{thm:obtuse-right}, in fact, follows from Theorem 1.1,
Theorem 13.4, and Corollary 13.5 in \cite{Annals} 

As in \S 7 of \cite{Annals},  for each 
critical point $p$ in a side $e$ of a triangle,
we choose coordinates so that $p=0$ and $e$ lies in the real-axis. 
In these coordinates, $u$ has a series expansion of the form 
$u(x+iy)= \sum_{jk} a_{jk} \cdot x^j \cdot y^k$.

\begin{prop}[compare Proposition 7.6 \cite{Annals}] 
\label{prop:tangential-cusp}
Suppose that $p$ is a degenerate critical point of $u$ in a side $e$.
If $a_{20} =0$ and $a_{30} \ne 0$, then $u$ has at least two critical
points in the boundary of $T$. Moreover, at least one of these critical 
points is a degree 1 vertex of $\Zcal(L_e u)$. 
\end{prop}

\begin{proof}
This follows from the first paragraph of the ``proof'' of Lemma 7.6
in \cite{Annals}.  In particular, since $a_{20}=0$ and $a_{30} \neq 0$
it follows that $\Zcal(L_e u) \setminus e$ contains has at least 
two components and hence  $\Zcal(L_e u)$ has two degree 1 vertices
that lie in $\partial T \setminus e$ (using Lemma 3.3 \cite{Annals}).
One of these could be the 
vertex opposite to $e$, but another one is a critical point 
distinct from $p$.
\end{proof}

\begin{prop}[compare Theorem 7.7 \cite{Annals}] \label{prep:thm}
If $u$ has a degenerate critical point $p$ that lies in
the boundary, then there exist at least two boundary critical points. 
Moreover, at least one critical point of $u$ in the boundary 
is a degree 1 vertex of the nodal set of some
directional or angular derivative of $u$. 
\end{prop}

\begin{proof}
The proof consists of combining the 
statement of Proposition \ref{prop:tangential-cusp}
with the statements of Propositions 7.3 and 7.4 \cite{Annals}.
\end{proof}

The  proof of Theorem \ref{thm:obtuse-right} is 
essentially a combination of the proofs of Theorem 1.1, Theorem 13.4, and, Corollary 13.5 in \cite{Annals}. 
However, to prove Theorem \ref{thm:obtuse-right} 
for right triangles, we will deviate from 
the exposition in \cite{Annals} a bit.
In \cite{Annals} we used 
the result in the acute case to prove 
the result for right triangles, and here 
we will first prove the result for obtuse 
triangles and then deduce the result 
for right triangles.

Let $t \mapsto u_t$ be the one-parameter family of eigenfunctions defined
at the beginning of \S9 \cite{Annals}. In particular, $u_0$
is a second Neumann eigenfunction on the given obtuse triangle $T_0$, 
and $u_1$ is a second Neumann eigenfunction on the right iscoceles triangle
with vertices $0$, $1$, and $i$.

As in \cite{Annals}, for each $t \in [0,1]$, 
let $N(t)$ denote the number of critical points of $u_t$.
The strategy is to show that Lemma 11.2 \cite{Annals}  
implies that if $u_0$ has at least two critical points, 
then $u_t$ has at least one critical point for each $t <1$ and 
sufficiently close to $1$.
\begin{proof}[Proof of Theorem \ref{thm:obtuse-right}]
If $u_0$ has an interior critical point, then Lemma 7.1 \cite{Annals} implies that $N(0) \geq 2$. 
Observe that the second Neumann eigenvalue 
of the right isoceles triangle
is simple, and so $u_1$ is a multiple of the 
function $u = \cos(\pi x) - \cos(\pi y)$. This implies that for $t < 1$ 
and sufficiently close to $1$, the eigenvalue 
$\mu_t$ is simple and $c_1(t) \neq 0$. Hence by the last part of Lemma 11.2 \cite{Annals} there exists a sequence $t_n$ converging to 
1 such that $N(t_n) \geq 2$ for each $n$.

Let $p$ be an accumulation point of a sequence of critical points, $p_n$, of $u_{t_n}$. 
If $p$ is not a vertex of the triangle $T_1$, then 
$p$ is a critical point of $u_1$. But $u_1$ is a multiple of the function $u = \cos(\pi x) - \cos(\pi y)$ that has no critical points. 
Therefore, each accumulation point of $\{p_n\}$ is a vertex. It follows from Proposition 9.1 \cite{Annals} that there exists $K>0$ such that
if $n >K$, then each critical point of $u_{t_n}$ lies in a side of the triangle.

Since $N(t_n) \geq 2$, Proposition 7.10 \cite{Annals}
implies that that for each $n>K$ there exist distinct sides $e$ and $e'$ 
and sequences of critical points $p_n$ and $p_n'$ so that for each
$n$ we have $p_n$ in  $e$ and $p_n'$ in $e'$. By passing to a subsequence 
if necessary we may assume that $p_n$ converges to a vertex 
$v$ of $T_1$, and $p_n'$ converges to a vertex $v'$ of $T_1$. By Lemma 9.2 \cite{Annals}, 
we have $v \neq v'$. The sets $\{v,v'\}$ and $\{1,i\}$ are both contained in 
a three element set, and hence we may assume without loss of generality that 
$v=1$ or $v=i$. Thus, by Lemma 9.3 \cite{Annals}, 
the function $u_1$ vanishes at either $1$ or $i$. 
But $u_1$ is a multiple of the 
function $u = \cos(\pi x) - \cos(\pi y)$ that does not vanish 
at $1$ or $i$.

Thus, each of the following
is true:
\begin{enumerate}
    \item $u_0$ has at most one critical point, and if such a critical point exists, then it lies in $\partial T_0$, and

    \item if $u_0$ has exactly one critical point then the first Bessel coefficient $c_1$ at the obtuse vertex $O$ of $T_0$ is non-zero.
\end{enumerate}
Finally we show that $u_0$ has no critical points. 
If $u_0$ had a critical point then, by above, $u_0$ has exactly one critical point and $c_1 \neq 0$. 
By Lemma 4.3 \cite{Annals} the obtuse vertex is not a local extrema of $u_0$.
But, by Proposition 8.1 \cite{Annals}, since $u_0$ has exactly one critical point each vertex of $T_0$ is a local extrema of $u_0$. 
This is a contradiction. 

This proves that if $T_0$ is an obtuse triangle then 
a second Neumann eigenfunction for $T_0$ has no critical points.
We next use this fact to prove the claim for right triangles. 

Given a (labeled) right triangle $T$, let $t \mapsto T_t$ 
be a path of labeled triangles 
such that $T_0=T$, and if $t >0$, then $T_t$ is an obtuse triangle. 
Let $u$ be a  $\mu_2(T_0)$-eigenfunction.
For each $t$, the eigenvalue  
$\mu_2(T_t)$ is simple \cite{Siudeja}, and hence 
standard perturbation theory implies that there exists a 
continuous path $t \mapsto u_t$ of $\mu_2(T_t)$-eigenfunctions,
such that $u_0=u$. If $u$ were to have an interior critical 
point, then a variant of Proposition 10.4 \cite{Annals}
would imply that $u_t$ has at least three critical points for small $t$. 
If $u$ were to have at least one degenerate critical point that belonged to a side of $T$, 
then a variant of Proposition \ref{prep:thm}
would imply that $u_t$ has
at least one critical points for $t$ small. 
Each is a contradiction  since $T_t$ is obtuse for $t > 0$. 
\end{proof}


\end{document}